\documentclass[a4paper,12pt,
reqno]{amsart}
\usepackage{latexsym,amscd,amssymb,amsmath,amsthm,comment,url}
\usepackage[dvips]{graphicx}
\usepackage[all]{xy}
\xyoption{poly}
\xyoption{knot}
\xyoption{import}
\makeatletter
\theoremstyle{plain}
  \newtheorem{thm}{Theorem}[section]
  \newtheorem{prop}[thm]{Proposition}
  \newtheorem{lem}[thm]{Lemma}
  \newtheorem{cor}[thm]{Corollary}
  
\theoremstyle{definition}
  \newtheorem{dfn}[thm]{Definition}
  \newtheorem{exmp}[thm]{Example}
\theoremstyle{remark}
  \newtheorem{rem}[thm]{Remark}
\let\opn\operatorname 
\let\term\emph
\def\@bothmode#1{\ifmmode #1\else $#1$\fi}
\newcount\@tempchn 
\def\@chCount#1{%
   \@tempchn=0
   \@tfor\member:=#1\do{\advance\@tempchn by 1}%
}
\def\@autopr#1{%
   \@chCount{#1}%
   \ifnum\@tempchn<2 #1\else (#1)\fi
}
\let\@tempopn\relax 
\def\@opform_#1#2{\@tempopn_{#1}\@autopr{#2}}
\numberwithin{equation}{section}
\def\NN{\mathbb{N}} 
\def\ZZ{\mathbb{Z}} 
\def\RR{\mathbb{R}} 
\def\kk{\Bbbk} 
\def\m{\ideal{m}} 
\def\p{\ideal{p}} 
\let\s\sigma 
\let\C\Sigma 
\let\t\tau 
\let\u\upsilon 
\def\MM{\mathcal M} 
\def\M{\mathbb M} 
\let\@tempar\relax 
\def\@seton^#1{\overset{#1}{\@tempar}}
\def\defar#1#2{\@xp\def\csname #1\endcsname{\def\@tempar{#2}\@ifnextchar^{\@seton}{\@tempar}}}
\defar{longto}{\longrightarrow} 
\defar{epito}{\twoheadrightarrow} 
\defar{monoto}{\rightarrowtail} 
\defar{embto}{\hookrightarrow} 
\def\imply{\@bothmode{\Rightarrow}} 
\def\Imply{\@bothmode{\Longrightarrow}} 
\def\iff{\@bothmode\Longleftrightarrow} 
\def\get{\@bothmode{\Leftarrow}} 
\def\Get{\@bothmode{\Longleftarrow}} 
\let\Dsum\bigoplus 
\def\sM{|\MM |} 
\def\supp{\opn{supp}} 
\def\cell{\mathcal X} 
\def\op{\mathsf{op}} 
\def\chara{\operatorname{char}} 
\def\P{\mathcal P}
\def\C{\mathcal C}
\def\om{\mathcal W}
\def\defopn#1{%
    \@xp\def\csname #1\endcsname{%
        \def\@tempopn{\opn{\csname the#1\endcsname}}%
        \@ifnextchar_{\@opform}{\@opform_{}}%
    }%
}

\defopn{Ker} 

\defopn{Im} 

\defopn{Cok} 

\defopn{Ann} 

\defopn{Ass} 

\defopn{Spec} 

\defopn{pd} 

\defopn{Sym} 

\defopn{id} 
\def\the@init{in}
\defopn{@init}
\def\init{\@init_{\succ}}
\def\relint{\opn{rel-int}} 
\let\ideal\mathfrak 
\def\theE{E}
\def\E{\@ifstar{{}^*\theE}{\theE}} 
\let\defcat\defopn
\def\theMod{Mod}    \def\themod{mod}
\def\theMod{Mod}    \def\themod{mod}

\let\the@Lgr\theMod  \let\the@lgr\themod
\defcat{Mod} 
\defcat{mod} 
\defcat{Gr} 
\defcat{gr} 
\defcat{@Lgr}
\def\Lgr#1{\@Lgr_{\ZZ\MM}{#1}} 
\defcat{@lgr}
\def\lgr#1{\@lgr_{\ZZ\MM}{#1}} 
\defcat{Sq} 
\defcat{InjSq} 
\def\theHom{Hom}    \def\theRHom{RHom}
\def\theExt{Ext}    \def\theTor{Tor}
\def\theD{D}
\let\@tempgrop\underline
\def\Hom{\@ifstar{\opn{\@tempgrop\theHom}}{\opn\theHom}} 
\def\RHom{\@ifstar{\opn{R\@tempgrop\theHom}}{\opn\theRHom}} 
\def\Ext{\@ifstar{\opn{\@tempgrop\theExt}}{\opn\theExt}} 

\def\Tor{\opn\theTor} 
\def\theDcat{{\mathsf D}}
\def\Db{\theDcat^b} 
\def\@G_#1{\Gamma_{#1}}
\def\G{\@ifnextchar_{\@G}{\@G_\m}} 
\def\DDD{\mathbf D}

\def\depth{\operatorname{depth}}
\def\rank{\operatorname{rank}}

\def\too{\longrightarrow}
\def\cpx#1{#1^{\bullet}} 
\def\theD{D} 
\def\D{\@ifstar{{}^*\theD}{\theD}} 
\def\cDx{\cpx {\mathcal D}_X} 
\def\const{\underline{\kk} } 

\def\cP{\mathcal P}
\def\cH{\mathcal H} 
\def\<{{\langle}}
\def\>{{\rangle}}

\def\ba{\mathbf a}
\def\bb{\mathbf b}

\def\b0{\mathbf 0}
\def\11{\mathbf 1}
\def\bA{\mathbf A}
\def\bH{\mathbf H}

\makeatother
\title[Higher Cohen-Macaulay property]
{Higher Cohen-Macaulay property of \\ squarefree modules and simplicial posets}

\author{Kohji Yanagawa}
\thanks{Partially supported by Grant-in-Aid for Scientific Research (c) (no.19540028).}
\address{Department of Mathematics, Kansai University,
Suita 564-8680, Japan}
\email{yanagawa@ipcku.kansai-u.ac.jp}
\begin{document}
%
\maketitle

\begin{abstract}
Recently, G. Fl\o ystad studied {\it higher Cohen-Macaulay property} 
of certain finite regular cell complexes.    
In this paper, we partially extend his results to 
squarefree modules, toric face rings, and simplicial posets.  
For example, we show that if (the corresponding cell complex of) 
a simplicial poset is $l$-Cohen-Macaulay then its codimension one skeleton 
is $(l+1)$-Cohen-Macaulay.  
\end{abstract}

\section{Introduction}
Let $\cell$ be a \term{finite regular cell complex} (cf. \cite[Section 6.2]{BH}). 
We use the convention that $\emptyset \in \cell$. 
Here the adjective ``regular" means that the closure $\overline{\s} \subset X$  
of an $i$-cell $\emptyset \ne \s \in \cell$ 
is always homeomorphic to the closed ball $\{ x \in \RR^i \mid ||x|| \leq 1 \}$. 
If we regard a finite simplicial complex as a cell complex, then it is regular. 

We say $\cell$ satisfies {\it the intersection property}, if  
for each $\s,\t \in \cell$, there is a cell $\u \in \cell$ such that 
$\overline \u = \overline \s \cap \overline \t$ (here, $\u$ can be $\emptyset$). 
Simplicial complexes and the boundary complexes of convex polytopes satisfy 
this property.

Clearly, the underlying space $X$ of $\cell$ has a finite simplicial decomposition. 
We say $\cell$ is Cohen-Macaulay over a field $\kk$ if the Stanley-Reisner ring 
$\kk[\Delta]$ is Cohen-Macaulay for some (equivalently, all) 
finite simplicial decomposition $\Delta$ of $X$. 
The Cohen-Macaulay property of $\kk[\Delta]$ may depend on $\chara(\kk)$. 
However, since we fix the base field $\kk$ throughout this paper, 
we often omit the phrase ``over $\kk$".   
We can define the Buchsbaum property of $\cell$ in the same way.
 
Let $V:= \{ \, \s \in \cell \mid \dim \s =0 \, \}$ be the set of vertices of $\cell$. 
For a subset $W \subset V$, set  
$$\cell |_W := \{ \, \t \in \cell \mid \s \subset \overline{\t}, \, \s \in V \Longrightarrow \s \in W \}.$$ 
This is a subcomplex of $\cell$. We simply denote $\cell|_{V \setminus W}$ by $\cell|_{-W}$. 

The following notion was introduced by Baclawski  (\cite{Bac}) in the case of simplicial complexes, 
and extended to regular cell complexes with the intersection property by 
Fl\o ystad (\cite{Fl}).

\begin{dfn}[c.f. \cite{Bac,Fl}]\label{l-CM def} 
For a positive integer $l$, 
we say $\cell$ is {\it $l$-Cohen-Macaulay} if $\cell |_{-W}$ is Cohen-Macaulay 
and $\dim \cell = \dim \cell |_{-W}$ for all $W \subset V$ with $\#W < l$.  
(Hence $\cell$ is 1-Cohen-Macaulay if and only if it is Cohen-Macaulay.)
\end{dfn}

Recently, Fl\o ystad gave the following striking results.

\begin{thm}[Fl\o ystad \cite{Fl}]\label{Flo}
Let $\cell$ be a finite regular cell complex with the intersection property. 
\begin{itemize}
\item[(1)] 2-Cohen-Macaulay property of $\cell$ is a topological property of 
the underlying space. 
\item[(2)] If $\cell$ is $l$-Cohen-Macaulay, then the codimension 1 skeleton 
$\cell' := \{ \, \s \mid \dim \s \leq \dim \cell -1 \, \}$  is $(l+1)$-Cohen-Macaulay.
\end{itemize}
\end{thm}

Motivated by the above result, we study the higher Cohen-Macaulay property 
(especially, 2-Cohen-Macaulay property) of relatively new notions 
of Combinatorial Commutative Algebra, such as squarefree modules, toric face rings, 
and (the face rings of) simplicial posets.  

For example, in Theorem~\ref{main for poset}, 
we show that Theorem~\ref{Flo} (2) holds for the corresponding 
regular cell complex $\Gamma(P)$ of a simplicial poset $P$, 
while $\Gamma(P)$ does not satisfy the intersection property 
and Theorem~\ref{Flo} (1) is no longer true.  

The notion of {\it toric face rings}, which generalizes both Stanley-Reisner rings 
and affine semigroup rings, is studied in (for example) \cite{BG, BKR, IR}. 
A toric face ring is supported by a finite regular cell complex 
satisfying the intersection property.  
In Theorem~\ref{2-CM}, we show that,   
under the assumption that a toric face ring $R$ is ``cone-wise normal", 
the supporting cell complex $\cell$ of $R$ is 2-Cohen-Macaulay 
if and only if $\cell$ (equivalently $R$) is Cohen-Macaulay and the 
canonical module $\omega_R$ of $R$ is generated by its degree 0 part.  
The corresponding statement does not hold for simplicial posets 
(at least, certain modification is required).

\section{Higher Cohen-Macaulay property of squarefree modules}
Let $S=\kk[x_1, \ldots, x_n]$ be a polynomial ring, and regard it as a $\ZZ^n$-graded ring. 
For $\ba = (a_1, \ldots, a_n) \in \NN^n$,
set $\supp(\ba) := \{ \, i \mid a_i \ne 0 \, \} \subset 
[n]:=\{1, \ldots ,n\}$,  and 
$x^\ba :=\prod x_i^{a_i} \in S$. 
Let $\gr S$ be the category of $\ZZ^n$-graded finitely generated $S$-modules and 
degree preserving $S$-homomorphisms.   
For $M \in \gr S$ and $\ba \in \ZZ^n$, 
$M_\ba$ denotes the degree $\ba$ component of $M$, and 
$M(\ba)$ denotes the shifted module of $M$ with $M(\ba)_\bb = M_{\ba + \bb}$. 
For $M, N \in \gr S$, $\Hom_S(M,N)$ is a 
$\ZZ^n$-graded $S$-module with $[\Hom_S(M,N)]_\ba = \Hom_{\gr S}(M, N(\ba))$.   
Hence $\Ext_S^i(M,N)$ has a similar $\ZZ^n$-grading. 

\begin{dfn}[\cite{Y}]
A $\ZZ^n$-graded $S$-module $M$ is called 
{\it squarefree}, if  it is finitely generated, $\NN^n$-graded 
(i.e., $M = \bigoplus_{\ba \in \NN^n} M_\ba$), and the multiplication map 
$M_\ba \ni y \longmapsto x^\bb y 
\in M_{\ba + \bb}$ is bijective for all $\ba, \bb \in \NN^n$ with 
$\supp(\ba) \supset \supp(\bb)$. 

Let $\Sq S$ be the full subcategory of $\gr S$ consisting of squarefree $S$-modules. 
\end{dfn}

For a simplicial complex $\Delta$ with the vertex set $[n]$, 
the Stanley-Reisner ring $\kk[\Delta] = S/I_\Delta$ is a squarefree $S$-modules.  
As shown in \cite{Fl,Y,Y04,Y09}, the notion of squarefree modules is 
useful in the study of the Stanley-Reisner rings.  
We have that $\Sq S$ is an abelian category with enough projectives 
and injectives, and an indecomposable injective is 
of the form $\kk[F] := S/(x_i \mid i \not \in F)$ for $F \subset [n]$.   

We say $\ba \in \NN^n$ is {\it squarefree} if  $a_i =0,1$ for all $i$. 
We sometimes identify a squarefree vector $\ba \in \NN^n$ with its support 
$\supp(\ba) \subset [n]$.  For example, in this situation, 
$M_{\supp(\ba)}$ denotes the homogeneous component $M_\ba$ of $M \in \gr S$.  
If $M \in \Sq S$, the  essential information of $M$ appears in 
its squarefree part $\bigoplus_{F \subset [n]} M_F$. 

For $W \subset [n]$ and $M \in \Sq S$, we can regard the $\kk$-vector space 
$$M|_W :=\bigoplus_{\substack{\ba \in \NN^n \\ \supp(\ba) \subset W}}M_\ba$$ 
as a squarefree module over the polynomial ring $S|_W := \kk[x_i \mid i \in W]$. 
This construction is a special case of \cite[Definition~3.26]{Mil}.   
We simply denote $M|_{[n] \setminus W}$ by $M|_{-W}$. 

Let $\kk[\Delta]$ be the Stanley-Reisner ring of a simplicial complex 
$\Delta$ with the vertex set $[n]$. Then, for $W \subset [n]$, we have 
$\kk[\Delta]|_W \cong \kk[\Delta|_W]$ as modules over $S|_W$, where 
$\Delta|_W:= \{ \, F \in \Delta \mid F \subset W \, \}$ 
is a simplicial complex with the vertex set $W$.

\begin{dfn}\label{tmk}
Let $l$ be a positive integer. We say $M \in \Sq S$ is {\it $l$-Cohen-Macaulay}, 
if for each $W \subset [n]$ with $\# W < l$, $M|_{-W}$  
is either the 0 module or a Cohen-Macaulay module with $\dim M|_{-W} = \dim M$. 
\end{dfn}

\begin{rem}\label{tmk rem}
(1) For $F \subset [n]$, set $\omega_F :=\kk[F](-F)$. Then $\omega_F$ is a squarefree module 
such that $(\omega_F)_{F'} \ne 0$ implies $F'=F$ 
(hence if $F \cap W \ne \emptyset$, then $\omega_F|_{-W} =0$). 
We see that $\omega_F$ is $l$-Cohen-Macaulay for all $l$.

(2) If $M \in \Sq S$ satisfies $M_0 \ne 0$, then $M|_{-W} \ne 0$ for all $W \subset [n]$. 
Hence the Stanley-Reisner ring $\kk[\Delta]$ is $l$-Cohen-Macaulay 
in the sense of Definition~\ref{tmk} if and only if $\Delta$ is $l$-Cohen-Macaulay 
in the usual sense (i.e., Definition~\ref{l-CM def}).  
\end{rem}

Let $M \in \gr S$. 
For $i \in \NN$ and $\ba \in \ZZ^n$, set $\beta^S_{i,\ba}(M) :=\dim_\kk [\Tor_i^S(\kk, M)]_\ba$. 
Then the $\ZZ^n$-graded minimal free resolution of $M$ is of the form 
$$\cdots \too \bigoplus_{\ba \in \ZZ^n}S(-\ba)^{\beta^S_{1,\ba}(M)} \too
\bigoplus_{\ba \in \ZZ^n}S(-\ba)^{\beta^S_{0,\ba}(M)} \too M \too 0.$$

\begin{lem}~\label{betti}
For $F \subset W \subset [n]$ and a squarefree module $M$, we have 
$$\beta_{i, F}^{S|_W}(M|_W) = \beta_{i, F}^S(M).$$ 
\end{lem}

\begin{proof}
By \cite[Proposition~3.27]{Mil}.  
\end{proof}

Let $\Db(\Sq S)$ be the bounded derived category of $\Sq S$. 
Set $\11 := (1,1, \ldots, 1) \in \NN^n$. Then $\omega_S:= S(-\11)$ is the 
$\ZZ^n$-graded canonical module of $S$, and its translation $\omega_S[n] \in \Db(\Sq S)$ 
gives a normalized $\ZZ^n$-graded dualizing complex of $S$. 
If $M$ is squarefree, then so is $\Ext_S^i(M, \omega_S)$.  
Hence $\DDD(-):=\RHom_S(-, \omega_S[n])$ gives a functor $\Db(\Sq S) \to \Db(\Sq S)^\op$.  
Note that $H^{-i}(\DDD(M))$ is isomorphic to $\Ext_S^{n-i}(M, \omega_S)$. 
As shown in \cite{Y04}, $\DDD(M) \in \Db(\Sq S)$ is quasi-isomorphic to the complex of the form 
\begin{equation}\label{DDD}
\cdots  \to \bigoplus_{\substack{F \subset [n] \\ \# F = i}} (M_F)^* \otimes_\kk \kk[F]
\to \bigoplus_{\substack{F \subset [n] \\ \# F = i-1}} (M_F)^* \otimes_\kk \kk[F]
\to \cdots \to (M_0)^* \otimes_\kk \kk  \to 0,
\end{equation}  
where $(M_F)^*$ denotes the dual $\kk$-vector space of $M_F$. 

If $\kk[\Delta]$ is a Buchsbaum (especially, Cohen-Macaulay) ring with 
$d:= \dim \kk[\Delta]$, then $\omega_{\kk[\Delta]} :=\Ext_S^{n - d}(\kk[\Delta], \omega_S)$ 
is called the {\it canonical module} of $\kk[\Delta]$.  This is a squarefree module. 
In the following sections, we will study a few generalizations of Stanley-Reisner 
rings. Their canonical modules are defined in a similar way.     

As shown in \cite{Mil,R0}, 
we can also define the {\it Alexander duality functor} $\bA : \Sq S \to (\Sq S)^\op$ 
as follows: For $M \in \Sq S$, $\bA(M)_F$ is the dual $\kk$-vector space of  
$M_{[n] \setminus F}$, and the multiplication map  
$\bA(M)_F \ni y \longmapsto x_i y \in \bA(M)_{F \cup \{ i \}}$ for $i \not \in F$ 
is the $\kk$-dual of $M_{[n] \setminus (F   \cup \{ i \})} 
\ni z \longmapsto x_i z \in M_{[n] \setminus F}$.   

\medskip

The next result is just a module version of \cite[Theorem~2.5]{Fl}. 
However, it has an application to (the face ring of) a simplicial poset.

\begin{thm}[{c.f. \cite[Theorem~2.5]{Fl}}]
\label{l-CM}
Let $M$ be a Cohen-Macaulay squarefree $S$-module with $\dim M = d$. 
For $l \geq 2$,  the following are equivalent. 
\begin{itemize}
\item[(i)] $M$ is $l$-Cohen-Macaulay. 
\item[(ii)] $\beta_{i,F}^S(M)=0$ for all $F \subset [n]$ and $i \in \NN$ such that 
$i > n -d-l+1$ and $\# F < i+d$. 
\item[(iii)] $\beta_{i,F}^S(\Ext_S^{n-d}(M,\omega_S))=0$ for all $F \subset [n]$ and $i \in \NN$ such that 
$i < l-1$ and $\# F > i$. 
\item[(iv)] $\bA (\Ext_S^{n-d}(M,\omega_S))$ satisfies Serre's condition $(S_{l-1})$.
\item[(v)] $\bA (\Ext_S^{n-d}(M,\omega_S))$ is the $(l-1)$-st syzygy module of some 
$S$-module.   
\end{itemize}
\end{thm}

\begin{proof} 
(i) $\Leftrightarrow$ (ii):
Since $M|_W$ is a squarefree module over $S|_W$,  $\beta_{i, \ba}^{S|_W}(M|_W) \ne 0$ 
implies that $\ba$ is a squarefree vector.  Hence, $M|_W$ is either the 0 module or 
a Cohen-Macaulay module of dimension $d$ 
if and only if $\beta^{S|_W}_{i,F}(M|_W)=0$ for all $F \subset W$ and 
$i \in \NN$ with $i > \#W -d$.  Hence the assertion follows from Lemma~\ref{betti} 
and easy computation.  

(ii) $\Leftrightarrow$ (iii): 
Let $P_\bullet$ be a $\ZZ^n$-graded minimal $S$-free resolution of $M$. 
Since $\omega_S \cong S(-\11)$ and $M$ is Cohen-Macaulay, $\Hom_S(P_\bullet, \omega_S)$ 
is a $\ZZ^n$-graded minimal free resolution of $\Ext_S^{n-d}(M, \omega_S)$ 
after suitable translation. Hence we have 
$$\beta_{i,F}^S(M) = \beta_{n-d-i,[n] \setminus F}^S(\Ext_S^{n-d}(M, \omega_S)).$$
Now the equivalence is clear. 

Note that if an $S$-module $N$ satisfies Serre's condition $(S_i)$ for $i \geq 1$ 
then $\Ass(N) = \Ass(S) =\{ \, (0) \, \}$ and $\dim N = n$.  
Hence we can prove the equivalence (iii) $\Leftrightarrow$ (iv) by the same way 
as \cite[Corollary~3.7]{Y}. 

The equivalence (iv) $\Leftrightarrow$ (v) is nothing other than 
(a special case of) \cite[Theorem~3.8]{EG}, which is a classical result 
essentially due to Auslander and Bridger. 
\end{proof}

The following is a special case of the equivalence (i) $\Leftrightarrow$ (iii) 
of Theorem~\ref{l-CM}. However, we remark it here, since this fact will be 
mentioned repeatedly in the following sections.   

\begin{cor}\label{2-CM sqf}
Let $M$ be a Cohen-Macaulay squarefree $S$-module with $\dim M=d$. 
Then $M$ is 2-Cohen-Macaulay if and only if $\Ext_S^{n-d}(M, \omega_S)$ 
is generated by its degree 0 part. 
\end{cor}

For $M \in \Sq S$ and $i \in \ZZ$ with $0 \leq i < d:= \dim M$, the submodule 
$$M^{>i}:= \bigoplus_{\substack{\ba \in \NN^n \\ \# \supp(\ba) > i}}M_\ba$$ 
of $M$ is a squarefree module again. Set $M^{<i>}:=M/M^{>i}$. 
Clearly,  $M^{\< i\>} \in \Sq S$ and $\dim M^{\<i\>} \leq i$ 
($M^{\< i \>}$ can be the 0 module). 
For a simplicial complex $\Delta$ with the vertex set $[n]$, 
let $\Delta^{(i)} = \{ F \in \Delta \mid \# F \leq  i+1 \}$ be its $i$-skeleton. 
Then, as an $S$-module, we have $\kk[\Delta^{(i)}] \cong \kk[\Delta]^{\< i+1 \>}$ 
(note that  $\dim \kk[\Delta] = \dim \Delta +1$).

\begin{thm}[{c.f. \cite[Corollary~2.7]{Fl}}]\label{main lemma}
Let $M$ be a squarefree $S$-module of dimension $d$, and $i$ an integer 
with $0 \leq i < d$. If $M$ is $l$-Cohen-Macaulay, then $M^{\< i \>}$ is 
$(l+d-i)$-Cohen-Macaulay (unless $M^{\< i\>} =0$).  
\end{thm}

\begin{proof}
Since $M$ is Cohen-Macaulay, only non-vanishing cohomology of the complex $\DDD(M)$ is 
$H^{-d}(\DDD(M))$, which is isomorphic to $\Ext_S^{n-d}(M, \omega_S)$. 
If we use \eqref{DDD} for a description of the complex $\cpx D:=\DDD(M)$, 
then $\DDD(M^{\< i\>})$ is the brutal truncation 
$\cdots \to 0 \to D^{-i} \to D^{-i+1} \to D^{-i+2} \to \cdots$
of $\cpx D$. Hence  
$H^j(\DDD(M^{\< i\>}))=0$ for all $j \ne -i$, that is, $M^{\< i \>}$ is Cohen-Macaulay.  
Since $H^{-i}(\DDD(M^{\< i\>}) \cong \Ext_S^{n-i}(M^{\< i\>}, \omega_S)$, 
we have the following exact sequences 
$$0 \to \Ext_S^{n-d}(M, \omega_S) \to D^{-d} \to D^{-d+1} \to \cdots \to 
D^{-i+1} \to  \Ext_S^{n-i}(M^{\< i\>}, \omega_S) \to 0$$
and 
$$0 \to \bA(\Ext_S^{n-i}(M^{\< i \>}, \omega_S)) \to \bA (D^{-i+1}) \to \cdots \to 
\bA(D^{-d}) \to  \bA(\Ext_S^{n-d}(M \omega_S)) \to 0.$$
Since $\bA(\kk[F])$ is isomorphic to $S(-([n]\setminus F))$, 
$\bA(D^j)$ is a free $S$-module for all $j$. Hence  $\bA(\Ext_S^{n-i}(M^{\< i \>}, \omega_S))$ 
is the $(d-i)$th syzygy of $\bA(\Ext_S^{n-d}(M, \omega_S))$. 
Hence the assertion follows from the equivalence between (i) and (v) of 
Theorem~\ref{l-CM}. 
\end{proof}

\section{Toric face rings and 2-Cohen-Macaulay cell complexes}
While we discuss {\it toric face rings} in this section, we only give ``casual" 
definition/construction of this ring and some related notions.  See \cite{OY} for precise information. 
The original construction found in \cite{BKR} is equivalent to that of \cite{OY}, 
but does not mention the regular cell complex  supporting a toric face ring.

A toric face ring is constructed from a {\it monoidal complex} $\MM$ supported by 
a finite regular cell complex $\cell$ with the intersection property. 
Here $\MM$ is a collection $\{ \M_\s \}_{\s \in \cell}$ of an affine semigroup  
$\M_\s \subset \ZZ^{\dim \s +1}$ (i.e., $\M_\s$ is a finitely generated additive submonoid 
of $\ZZ^{\dim \s+1}$) with $\ZZ\M_\s = \ZZ^{\dim \s +1}$ and 
$\M_\s \cap (-\M_\s) = \{ 0 \}$.  Of course, we require several conditions on $\M_\s$'s  
(not all  finite regular cell complexes with the intersection property can support a monoidal complex).  
We assume that the boundary complex of the cross section 
$\overline{\cP}_\s$ of the polyhedral cone $\cP_\s := \RR_{\geq 0} \M_\s \subset \RR^{\dim \s +1}$ 
can be identified with the subcomplex $\{ \t \mid \t \subset \overline{\s} \}$ of $\cell$
(note that $\overline{\cP}_\s$ is a convex polytope of dimension $\dim \s$). 
The face $\C_\t$ of $\cP_\s$ corresponding to $\t \subset \overline{\s}$ 
is isomorphic to $\cP_\t$ as a polyhedral cone. Moreover, the monoid 
$\M_\t$ is isomorphic to $\C_\t \cap \M_\s$.

Let $|\MM|$ be the set given by glueing all $\M_\s$'s along with $\cell$. 
We can regard $\M_\s \subset |\MM|$ for all $\s \in \cell$. For $a, b \in |\MM|$,  
we can not define their sum in general, 
that is, $|\MM|$ is no longer a monoid. 
However, if there is some $\s \in \cell$ with $a,b \in \M_\s$, 
then we have their sum $a+b \in \M_\s \subset |\MM|$.

Then the toric face ring  $R:= \kk[\MM]$ of $\MM$ over $\kk$ is the vector space 
$\Dsum_{a \in \sM} \kk \, t^a$ with the $\kk$-linear multiplication defined by  
$$
t^a \cdot t^b = \begin{cases}
t^{a+b} & \text{if $ a,b \in \M_\s $ for some $ \s \in \cell $;}\\
0       & \text{otherwise.}
\end{cases}
$$
Note that $\dim R = \dim \cell +1$. 

For each $\s \in \cell$, $\p_\s := ( \, t^a \mid 
a \not \in \M_\s \, )$ is a prime ideal of $R$ with $R/\p_\s \cong \kk[\M_\s]$, where 
$\kk[\M_\s]$ is the semigroup ring of $\M_\s$. 
Clearly, $\kk[\M_\s]$ can be seen as a subring of $R$.  
In $R$,  $\kk[\M_\s]$ for $\s \in \cell$ are ``glued" along with $\cell$.  
We say $\MM$ is {\it cone-wise normal} if $\kk[\M_\s]$ is normal 
(equivalently, $\ZZ^{\dim \s +1}\cap \RR_{\geq 0} \M_\s =\M_\s$) for all $\s \in \cell$.

\begin{exmp}\label{prototypical}
Let $A:= \kk[ \, x^\ba \mid \ba \in \M \, ] \subset \kk[ \, x_1, \ldots, x_n \,]$ 
be the semigroup ring of an affine semigroup $\M \subset \NN^n$. An easy example of a 
toric face ring is the quotient ring $A/I$ of $A$ by a radical $\ZZ^n$-graded ideal $I$. 
Let $\P:= \RR_{\geq 0} \M \subset \RR^n$ be the polyhedral cone spanned by $\M$, 
and $L$ its face lattice.  For the ideal $I$, 
there is a subset $\Sigma$ of $L$ such that $$A/I = \bigoplus_{\substack{\ba \in \C \cap \ZZ^n \\ 
\C \in \Sigma}} \kk \, x^\ba.$$
(Note that $\C \in \Sigma$ is a face of $\P$, and it is also a cone in $\RR^n$.)
Clearly, $\{ \C \}_{\C \in \Sigma}$ forms a polyhedral fan, that is, 
$\C' \subset \C \in \Sigma$ and $\C' \in L$ imply that $\C' \in \Sigma$.  
The monoidal complex giving $A/I$ is $\MM := \{ \, \M \cap \C \mid \C \in \Sigma  \, \}$. 
Let $H \subset \RR^n$ be a hyperplane intersecting $\P$ transversely. 
The cell complex supporting $\MM$ is given by $\{\,  \relint(H \cap \C) \mid \C \in \Sigma \, \}$. 
If $A$ is normal, then $A/I$ is cone-wise normal. 

If $A$ is the polynomial ring $\kk[x_1, \ldots, x_n ]$, then $A/I$ can be attained as 
the Stanley-Reisner ring $\kk[\Delta]$ of a simplicial complex $\Delta$. 
In this case, the supporting cell complex is nothing other than $\Delta$. 
\end{exmp}

Clearly, $A/I$ of the above Example has a  $\ZZ^n$-grading inherited from that of $A$. 
However, a toric face ring does not admit a nice multi-grading in its most general setting, 
while the decomposition $R=\Dsum_{a \in \sM} \kk \, t^a$ plays a similar role to the 
graded structure.  

\medskip

\noindent{\bf Known Results.} 
Let $\MM$ be a cone-wise normal monoidal complex supported by a cell complex $\cell$. 
For $R:= \kk[\MM]$, we have the following. 
See \cite{OY} for detail. 
\begin{itemize}
\item[(1)] 
$R$ is Cohen-Macaulay (resp. Buchsbaum) if and only if $\cell$ is Cohen-Macaulay 
(resp. Buchsbaum) in the sense of \S1.

\item[(2)] We can naturally define squarefree modules over $R$. 
For example, $R$ itself is squarefree. 
If $R$ is Buchsbaum, then the canonical module $\omega_R$ of $R$ is also. 

A squarefree $R$-module $M$ has the decomposition $M = \bigoplus_{a \in \sM} M_a$ 
as a $\kk$-vector space. Note that $\sM$ has the 0 element. In the sequel,   
the homogeneous component $M_0$ plays a role. 

\item[(3)]  A squarefree $R$-module $M$ gives the constructible sheaf $M^+$
on the underlying space $X$ of $\cell$.  
For example, $R^+$ is the $\kk$-constant sheaf $\const_X$ on $X$.   
\end{itemize}

Let $\MM$ be a cone-wise normal monoidal complex supported by a regular cell complex $\cell$ 
with $d=\dim \cell$. 
Then the underlying space $X$ of $\cell$ admits the Verdier dualizing complex $\cDx$ 
with the coefficients in $\kk$.  
By \cite[Theroem~6.4]{OY}, $R:=\kk[\MM]$ is Buchsbaum if and only if 
$\cH^i(\cDx)=0$ for all $i \ne -d$. 
In this case, set $$\om_X := \cH^{-d}(\cDx).$$ 
If $X$ is a manifold, then $R$ is Buchsbaum and $\om_X$ is 
the {\it orientation sheaf} of $X$ with the coefficients in $\kk$. 
Moreover, $(\omega_R)^+ \cong \om_X$.

\begin{lem}\label{global section}
Let $R$ be a cone-wise normal toric face ring, and  $M$ a squarefree module over  
$R$ with $\depth_R M \geq 2$. 
Then $M$ is generated by its degree 0 part $M_0$ if and only if the sheaf $M^+$ is generated 
by global sections.  
\end{lem}

\begin{proof}
For each $\emptyset \ne \s \in \cell$, take an element $a(\s) \in \M_\s \subset |\MM|$ 
contained in the interior of the cone $\RR_{\geq 0} \M_\s$ in $\RR^{\dim \s +1}$. 
The stalk $(M^+)_p$ at any point $p \in  \s$ is isomorphic to 
$M_{a(\s)}$. 
Since $\depth_R M \geq 2$, we have $\Gamma(X,M^+) \cong M_0$ by \cite[Theorem~6.2 (a)]{OY},  
and the natural map $\Gamma(X, M^+) \to (M^+)_p$ corresponds to the map 
$\varphi_\s: M_0 \ni x \longmapsto t^{a(\s)}x \in M_{a(\s)}$. Hence $M^+$ is generated by global sections 
if and only if the map $\varphi_\s$ is surjective for all $\s \ne \emptyset$. 
Since $M$ is squarefree, the latter condition states that 
$M$ is generated by $M_0$ as an $R$-module. 
\end{proof}

 
\begin{thm}\label{2-CM}
Le $\MM$ be a cone-wise normal monoidal complex supported by a cell complex $\cell$ 
with the underlying space $X$. 
For the toric face ring $R:= \kk[\MM]$, the following are equivalent.
\begin{itemize}
\item[(i)] $\cell$ is 2-Cohen-Macaulay. 
\item[(ii)] $R$ is Cohen-Macaulay, 
and the canonical module $\omega_R$ is generated by its degree 0 part. 
\end{itemize}

If $\dim R \geq 2$, the above conditions are also equivalent to the following. 
\begin{itemize}
\item[(iii)] $\cell$ is Cohen-Macaulay and the sheaf $\om_X$ is generated by  global sections. 
\end{itemize}
\end{thm}

\begin{proof}
A toric face rings of dimension 1 (i.e., the case when $\dim \cell =0$) 
is always a Stanley-Reisner ring. Moreover, in this case, $\cell$ is 2-Cohen-Macaulay 
unless $\cell$ consists of a single point.  Hence the assertion 
is easy if $\dim R =1$. 

So it suffices to show the equivalence of (i), (ii), and (iii) 
under the assumption that $\dim R \geq 2$. Now $\depth_R \omega_R \geq 2$. 
Since $\om_X \cong (\omega_R)^+$, the equivalence (ii) $\Leftrightarrow$ (iii) 
follows from Lemma~\ref{global section}. 
It remains to prove the equivalence (i) $\Leftrightarrow$ (iii). 
Let $\Delta$ be the barycentric subdivision of $\cell$. 
Since $\cell$ satisfies the intersection property, 
$\cell$ is 2-Cohen-Macaulay if and only if so is $\Delta$ by \cite[Theorem~2.8]{Fl}. 
On the other hand, as shown in \cite{Bac} 
(and as a special case of Corollary~\ref{2-CM sqf}), 
$\Delta$ is 2-Cohen-Macaulay if and only if $\Delta$ is Cohen-Macaulay and 
the canonical module $\omega_{\kk[\Delta]}$ of the Stanley-Reisner ring 
$\kk[\Delta]$ is generated by its degree 0 part. 
Since $\depth \omega_{\kk[\Delta]} \geq 2$ and $(\omega_{\kk[\Delta]})^+ \cong \om_X$, 
the desired equivalence follows from Lemma~\ref{global section}.  
\end{proof}

\begin{rem}
Let $\cell$ be a finite regular cell complex with the intersection property, 
and set $V:= \{ \, \s \in \cell \mid \dim \s = 0 \, \}$ and $d:= \dim \cell$. 
Fl\o ystad (\cite{Fl}) constructed the $i$th {\it enriched cohomology} 
$\bH^i(\cell; \kk)$ (or just $\bH^i(\cell)$, since we have fixed the base field $\kk$) 
of $\cell$, which is a squarefree module over the polynomial ring 
$S:= \kk[\, x_\s \mid \s \in V \, ]$. If $\Delta$ is a simplicial complex, 
we have $\bA(\bH^i(\Delta)) \cong \Ext_S^{\# V-i-1}(\kk[\Delta], \omega_S)$ 
for all $i$.  Even for general $\cell$, it is Cohen-Macaulay if and only if $\bH^i(\cell)=0$ 
for all $i \ne d$. 

From now, we assume that $\cell$ is Cohen-Macaulay. 
\cite[Theorem~2.4]{Fl} states that $\cell$ is $l$-Cohen-Macaulay 
if and only if $\bH^d(\cell)$ is the $(l-1)$-st syzygy module of some 
$S$-module.  Hence $\cell$ is 2-Cohen-Macaulay if and only if 
$\bA(\bH^d(\cell))$ is generated by its degree 0 part.  
Clearly, this is analogous to Theorem~\ref{2-CM}. 
However, the relation between $\bA(\bH^d(\cell))$ and $\omega_R$ is not direct 
unless $\cell$ is a simplicial complex.  
\end{rem}

\section{Simplicial Poset}
A finite partially ordered set ({\it poset}, for short) 
$P$ is called {\it simplicial}, if it admits the smallest element 
$\hat{0}$, and the interval $[\hat{0},  x]:= \{ \, y \in P \mid y \leq x \, \}$ is isomorphic to 
a boolean algebra for all $x \in P$. 
For the simplicity, we denote $\rank (x)$ of $x \in P$ by $\rho(x)$. 
If $P$ is simplicial and $\rho(x) = m$, then $[\hat{0}, x]$ is isomorphic to the boolean algebra 
$2^{\{1, \ldots, m\}}$. 

Let  $\Delta$ be a finite simplicial complex (with $\emptyset \in \Delta$). 
Its face poset (i.e., the set of the faces of $\Delta$ with the order given by inclusion) 
is a simplicial poset.   
Any simplicial poset $P$ is the face (cell) poset of a regular cell complex, 
which we denote by $\Gamma(P)$. Unless $\Gamma(P)$ is a simplicial complex, 
$P$ does not satisfy the intersection property.  
For example, if two $d$-simplices are glued along their boundaries, then it is not a simplicial 
complex, but gives a simplicial poset.  

From now on, let $P$ be a simplicial poset. For $x,y \in P$, set 
$$[x \vee y]:= \text{the set of minimal elements of $\{ \, z \in P \mid z \geq x, y \, \}$.}$$
More generally, for $x_1, \ldots, x_m \in P$, $[x_1 \vee \cdots \vee x_m ]$  
denotes the set of minimal elements of the common upper bounds of  $x_1,  \ldots, x_m$. 

Set $V:=\{ \, y \in P \mid \rho(y) = 1 \, \} = \{y_1, \ldots, y_n \}$. 
We sometimes identify $V$ with $[n]$ in the natural way. 
For $U \subset V$, we simply denote $[\bigvee_{i  \in U} y_i]$ by $[U]$. 
Then we have $P = \coprod_{U \subset V}[U].$
For $W \subset V$, set 
$$P|_W := \coprod_{U \subset W}[U].$$
Clearly, $P|_W$ is simplicial again. 
We simply denote $P|_{V \setminus W}$ by $P|_{-W}$. 

Stanley \cite{St91} defined the face ring $A_P$ of a simplicial poset $P$. 
For the definition, we remark that if $[ x \vee y] \ne \emptyset$ then  
$\{ \, z \in P \mid z \leq x, y \, \}$ has the largest element, which is denoted by $x \wedge y$.   
Let $S := \kk[ \, t_x \mid x \in P \, ]$ be 
the polynomial ring in the variables $t_x$. Consider the ideal  
$$I_P:= (\, t_xt_y - t_{x \wedge y} \sum_{z \in [x \vee y]} t_z \ | \ x,y \in P   \, ) 
+ (\, t_{\hat{0}} -1 \, )$$
of $S$ (if $[x \vee y] = \emptyset$, we interpret that 
$t_{x \wedge y} \cdot \sum_{z \in [x \vee y]} t_z = 0$), and set $$A_P:= S/I_P.$$ 
We have  $\dim A_P = \rank P = \dim \Gamma(P) +1$. 

For a rank 1 element $y_i \in V$, we simply denote $t_{y_i}$ by $t_i$.  
If $\{ x \}= [U]$ for some $U \subset V$ with $\# U \geq 2$, then $t_x = \prod_{i \in U} t_i$ 
in $A_P$, and $t_x$ is a ``dummy".  Since $I_P$ is a homogeneous ideal under the 
grading given by $\deg(t_x) = \rho(x)$, $A_P$ is a graded ring. The algebra 
$A_P$ is generated by degree 1 elements if and only if $\Gamma(P)$ is a simplicial complex. 
In this case, $A_P$ coincides with the Stanley-Reisner ring of $\Gamma(P)$.

We say $P$ is Cohen-Macaulay (resp. Buchsbaum), if  the cell complex 
$\Gamma(P)$ is  Cohen-Macaulay (resp. Buchsbaum) in the sense of \S1. 
Duval~\cite{Du} showed that 
$P$ is Cohen-Macaulay if and only if $A_P$ is a Cohen-Macaulay ring. 
The same is true for the Buchsbaum property (c.f. \cite{Y09}).  

\begin{prop}\label{tomo}
Some (equivalently, all) finite simplicial decomposition of $\Gamma(P)$ 
is 2-Cohen-Macaulay if and only if $A:=A_P$ is Cohen-Macaulay and the canonical module
$\omega_A$ is generated by its degree 0 part. 
\end{prop}

\begin{proof}
We can define squarefree modules over $A$ in a natural way, 
and a squarefree $A$-module $M$ gives the constructible sheaf $M^+$ on 
$\Gamma(P)$ (see \cite{Y09}).  
The theory of squarefree modules over $A$ is quite parallel to that over a toric face ring.  
Hence we can prove the assertion by the same way as Theorem~\ref{2-CM}. 
\end{proof}

\begin{dfn}
We say a simplicial poset $P$ is {\it $l$-Cohen-Macaulay} if $P|_{-W}$ is Cohen-Macaulay 
and $\rank (P|_{-W}) = \rank P$ for all $W \subset  V$ with $\# W < l$. 
Clearly, $P$ is $l$-Cohen-Macaulay if and only if the cell complex $\Gamma(P)$ is 
$l$-Cohen-Macaulay. 
\end{dfn}

Recall our convention that $V=\{ \, y \in P \mid \rho(y) = 1 \, \} = \{y_1, \ldots, y_n \}$.
Clearly, $A:=A_P$ has a $\ZZ^n$-grading such that $\deg t_i \in \NN^n$ is the $i$th unit vector. 
Consider the polynomial ring $T:= \Sym (A_1)  \cong \kk[t_1, \ldots, t_n]$.  
Then $A$ is a finitely generated $\ZZ^n$-graded $T$-module, 
moreover $A$ is a squarefree $T$-module.   
If $A$ is Buchsbaum, then $\omega_A$ is a squarefree $T$-module.   

Let $W \subset V$. 
If we regard the face ring $A_{P|_W}$ of $P|_W$ as a module over 
the polynomial ring $T|_W = \kk[t_i \mid i \in W]$, 
it coincides with the restriction $(A_P)|_W$ 
of $A_P$ as a squarefree $T$-module  (see \S2).  
As Remark~\ref{tmk rem}, $P$ is $l$-Cohen-Macaulay if and only if  
$A_P$ is $l$-Cohen-Macaulay as a squarefree $T$-module. 

\begin{prop}
Let the notation be as above. A simplicial poset 
$P$ is $2$-Cohen-Macaulay if and only if $A:=A_P$ is Cohen-Macaulay and 
the canonical module $\omega_A$ is generated by its degree 0 part as a $T$-module. 

In particular, if $P$ is 2-Cohen-Macaulay then some (equivalently, all) 
finite simplicial decomposition of $\Gamma(P)$ is 2-Cohen-Macaulay.  
\end{prop}

\begin{proof}
The first assertion follows from Corollary~\ref{2-CM sqf}. 
The second follows from the first and Proposition~\ref{tomo}. 
\end{proof}

\begin{rem}
Even if a finite simplicial decomposition of $\Gamma(P)$ is 2-Cohen-Macaulay, 
$P$ itself is not so in general. 
For example, let $P$ be the simplicial poset given by two $d$-simplices  
glued along their boundaries. Since the underlying space of 
$\Gamma(P)$ is a $d$-dimensional sphere, 
its simplicial decompositions are 2-Cohen-Macaulay. 
However $\rank (P|_{-\{y\}}) < \rank P$ for all $y \in V$, 
and $P$ is not 2-Cohen-Macaulay. 
\end{rem}

Consider the induced subposet 
$P^{\< i \>}:= \{ \, \, x \in P \mid \rho(x) \leq i \, \}$ 
of $P$. This is a simplicial poset again. 

\begin{thm}[{c.f. \cite[Corollary~2.7]{Fl}}]\label{main for poset} 
Let $P$ be a simplicial poset with $\rank P = d$. 
If $P$ is $l$-Cohen-Macaulay then $P^{\<i\>}$ is $(l+ d-i)$-Cohen-Macaulay. 
\end{thm}

\begin{proof}
Recall that $A_P$ is a squarefree module over 
$T = \kk[ \, t_1, \ldots, t_n \,]$. 
The face ring $A_{P^{\<i \>}}$ of $P^{\< i\>}$ 
coincides with the skeleton  $(A_P)^{\<i\>}$ of $A_P$ as a squarefree $T$-module. 
Hence the assertion 
follows from Theorem~\ref{main lemma}. 
\end{proof}

The next statement immediately follows from Theorem~\ref{main for poset}: 
Let $P$ be a Cohen-Macaulay simplicial poset  with $\rank P =d$. 
Then the edge graph of $P$ (i.e., the skeleton $P^{\<2\>}$) 
is $(d-1)$-connected. However, this also follows from \cite[Theorem~4.5]{Du}. 

\begin{rem}[Buchsbaum* complex] 
Athanasiadis and Welker (\cite{AW}) call a finite simplicial complex $\Delta$ 
{\it Buchsbaum*} (over $\kk$), if $\kk[\Delta]$ is Buchsbaum and the canonical module 
$\omega_{\kk[\Delta]}$ is generated by its degree 0 part. 
(Their original definition is different,  but equivalent to the above one by   
 \cite[Proposition~2.8]{AW}). Since the Buchsbaum* property is topological, 
we say a finite regular cell complex $\cell$ is Buchsbaum*, 
if some (equivalently, all) finite simplicial decomposition of $\cell$ is so.  
By the same argument as the proof of Theorem~\ref{2-CM}, we have the following. 

\begin{itemize}
\item[(1)] Let $R$ be a cone-wise normal toric face ring. 
The supporting cell complex of $R$ is  Buchsbaum* if and only if $R$ is 
Buchsbaum and the canonical module $\omega_R$ is generated by its degree 0 part. 
\item[(2)] Let $P$ be a simplicial poset, and $A$ its face ring. 
Then $\Gamma(P)$ is Buchsbaum* if and only if  
$A$ is Buchsbaum and $\omega_A$ is generated by its degree 0 part. 
\end{itemize}
\end{rem}

\bigskip

\section*{Acknowledgement} 
The author is grateful to Ezra Miller for his helpful comments on an earlier version 
of this paper.


\begin{thebibliography}{99}
\bibitem{AW} C.A. Athanasiadis and V. Welker, Buchsbaum* complexes, preprint ({\tt arXiv:0909.1931}). 

\bibitem{Bac} K. Baclawski,  Cohen-Macaulay connectivity and geometric lattices, 
 Europ. J. Combinatorics {\bf 3}  (1982), 293--305.

\bibitem{BG}
W. Bruns and J. Gubeladze, Polytopes, rings, and $K$-theory, Springer, 2009. 

\bibitem{BH} W. Bruns and J. Herzog,  ``Cohen-Macaulay rings", revised edition,  
Cambridge University Press, 1998.

\bibitem{BKR}
W. Bruns, R. Koch, and T. R\"omer,
\term{Gr\"obner bases and Betti numbers of monoidal complexes}, 
Michigan Math. J.  {\bf 57} (2008), 71-91. 

\bibitem{Du} A.M. Duval, Free resolutions of simplicial posets,  
J. Algebra {\bf 188} (1997), 363--399. 

\bibitem{EG} E. G. Evans and P. Griffith, 
Syzygies, London Mathematical Society Lecture Note Series, vol. 106, 1985.


\bibitem{Fl} G. Fl\o ystad, Cohen-Macaulay cell complexes, in Algebraic and Geometric
Combinatorics, C. A. Athanasiadis et al.  eds.,
Contemporary Mathematics vol. 423,  American Mathematical Society (2007), pp. 205--220. 


\bibitem{IR} B. Ichim and T. R\"omer,
\term{On toric face rings}, J. Pure Appl. Algebra {\bf 210} (2007), 249--266.

\bibitem{Mil} E. Miller, 
The Alexander duality functors and local duality with monomial support, 
J. Algebra. {\bf 231} (2000), 180--234. 


\bibitem{OY} R. Okazaki and K. Yanagawa, \term{Dualizing complex of a toric face ring,} 
Nagoya Math. J. {\bf 196} (2009), 87--116. 


\bibitem{R0} T. R\"omer,  Generalized Alexander duality and applications, 
Osaka J. Math. {\bf 38} (2001), 469--485.  


\bibitem{St91} 
R. Stanley, $f$-vectors and $h$-vectors of simplicial posets, 
J. Pure Appl. Algebra {\bf 71} (1991), 319--331. 


\bibitem{St} R. Stanley,  
``Combinatorics and commutative algebra", 2nd ed. Birkh\"auser 1996.   


\bibitem{Y} K. Yanagawa, 
Alexander duality for Stanley-Reisner rings and squarefree 
$\NN^n$-graded modules, J. Algebra 225 (2000), 630--645.  


\bibitem{Y04} K. Yanagawa, 
Derived category of squarefree modules and local cohomology with monomial 
ideal support, {\it J. Math. Soc. Japan} {\bf 56} (2004) 289--308. 


\bibitem{Y09}
K. Yanagawa, 
Dualizing complex of the face ring of a simplicial poset, preprint 
({\tt arXiv:0910.1498}). 
\end{thebibliography}
\end{document}